\documentclass[12pt,reqno,twoside]{amsart}

\usepackage{amsfonts}
\usepackage{amsmath,amssymb,amsthm,amsthm}
\usepackage{mathtools}
\usepackage{dsfont}
\usepackage{etoolbox}
\usepackage{color}
\usepackage{graphicx,color}
\usepackage{dsfont,subfigure}
\usepackage{fixmath}
\usepackage{tcolorbox}
\usepackage{algorithm}
\usepackage{algpseudocode}
\usepackage{graphicx}
\usepackage{enumerate}
\usepackage{tabularx}
\usepackage{listings}

\usepackage[dvips,margin=1.5in]{geometry}

\makeatletter
\def\eqnarray{\stepcounter{equation}\let\@currentlabel=\theequation
\global\@eqnswtrue
\tabskip\@centering\let\\=\@eqncr
$$\halign to \displaywidth\bgroup\hfil\global\@eqcnt\z@
  $\displaystyle\tabskip\z@{##}$&\global\@eqcnt\@ne
  \hfil$\displaystyle{{}##{}}$\hfil
  &\global\@eqcnt\tw@ $\displaystyle{##}$\hfil
  \tabskip\@centering&\llap{##}\tabskip\z@\cr}

\def\endeqnarray{\@@eqncr\egroup
      \global\advance\c@equation\m@ne$$\global\@ignoretrue}

\newtheorem{thm}{Theorem}[section]

\newtheorem{prop}[thm]{Proposition}

\newtheorem{remark}[thm]{Remark}

\usepackage[textsize=small]{todonotes}
\title[Deep Neural Nets With Fixed Bias Configuration]
{Deep Neural Nets with Fixed \\ Bias Configuration}

\author[H.~Antil, T.S. Brown, R.~L\"ohner, F.~Togashi, D.~Verma]
{Harbir Antil$^{1}$, Thomas S. Brown$^{1,2}$, Rainald L\"ohner$^{2}$, Fumiya Togashi$^{3}$, Deepanshu Verma$^{4}$}

\address{$^1$Center for Mathematics and Artificial Intelligence (CMAI),
             College of Science, 
             George Mason University,
             Fairfax, VA 22030-4444, USA.}
\email{hantil@gmu.edu,tbrown62@gmu.edu}

\address{$^2$Center for Computational Fluid Dynamics,
             College of Science, George Mason University, 
             Fairfax, VA 22030-4444, USA.}
\email{rlohner@gmu.edu}               

\address{$^3$Applied Simulations, Inc.,
             1211 Pine Hill Road, McLean, VA 22101, USA}
\email{fumiya.togashi@gmail.com}               
     
\address{$^4$Department of Mathematics, 
	       Emory University, Atlanta, GA 30322, USA}
\email{deepanshu.verma@emory.edu}

\thanks{
This work is partially supported by the Defense Threat
Reduction Agency (DTRA) under contract HDTRA1-15-1-0068 where 
Jacqueline Bell served as the technical monitor, and by NSF grants 
DMS-2110263, DMS-1913004, and the Air Force Office of Scientific 
Research under Award NO: FA9550-19-1-0036.}

\begin{document}

\begin{abstract}
For any given neural network architecture a permutation of weights 
and biases results in the same functional network. This implies
that optimization algorithms used to `train' or `learn' the network
are faced with a very large number (in the millions even for small
networks) of equivalent optimal solutions in the parameter space. To 
the best of our knowledge, this observation is absent in the literature.
In order to narrow down the parameter search space, a novel technique 
is introduced in order to fix the bias vector configurations to be 
monotonically increasing. This is achieved by augmenting a typical 
learning problem with inequality constraints on the bias vectors in each 
layer. A Moreau-Yosida regularization based algorithm is proposed to 
handle these inequality constraints and a theoretical convergence of 
this algorithm is established. Applications of the proposed approach 
to standard trigonometric functions and more challenging stiff ordinary 
differential equations arising in chemically reacting flows clearly 
illustrate the benefits of the proposed approach.  
Further application of the approach on the MNIST dataset within 
TensorFlow, illustrate that the presented approach can be incorporated 
in any of the existing machine learning libraries.

\end{abstract}

\maketitle

\section{Introduction}\label{s:intro}

\medskip
\noindent
{\bf Background.} 
A typical neural network can be represented 
as a function $\mathcal F: \mathbb R^{n_0} \to \mathbb R^{n_L}$ that consists 
of the composition of layer functions $\{f_\ell\}_{\ell = 0}^{L-1}$ and 
can written as 
\begin{equation} \label{eq:NNFunc}
	\mathcal F = f_{L-1} \circ f_{L-2} \circ \cdots \circ f_0 \, . 
\end{equation}
Each layer function is parameterized by a weight matrix 
$W_\ell \in \mathbb R^{n_\ell \times n_{\ell+1}}$, a bias vector 
$b_\ell \in \mathbb R^{n_{\ell+1}}$, and incorporates a nonlinear activation 
function $\sigma$, for instance, ReLU \cite{MR3617773}.

The weights $W_\ell$ and biases $b_\ell$ are determined in a process known 
as training the network. 
Training a neural network can be written in the framework of constrained 
optimization as follows: for training data $\{u^i, S(u^i)\}_{i=1}^N$ 
(input/output pairs), solve 
	\begin{equation}\label{eq:orgprob}
	\begin{aligned}
		&\min_{\{W_\ell\}_{\ell=0}^{L-1}, \{b_\ell\}_{\ell=0}^{L-2}} J(\{(y_L^i, S(u^i))\}_i,\{W_\ell\}_\ell,\{b_\ell\}_\ell) 
		\\
		 &\mbox{subject to } \quad y_L^i = \mathcal F(u^i; (\{W_\ell\}, \{b_\ell\})) \qquad i = 1, \dots, N,  
	\end{aligned}
	\end{equation}
where the function $J$, known as the loss function, measures the error of 
the approximation of $S(u^i)$ by the network output $y_L^i$ in some way.  
 
\begin{figure}[h!]
\centering
\includegraphics[width = \textwidth]{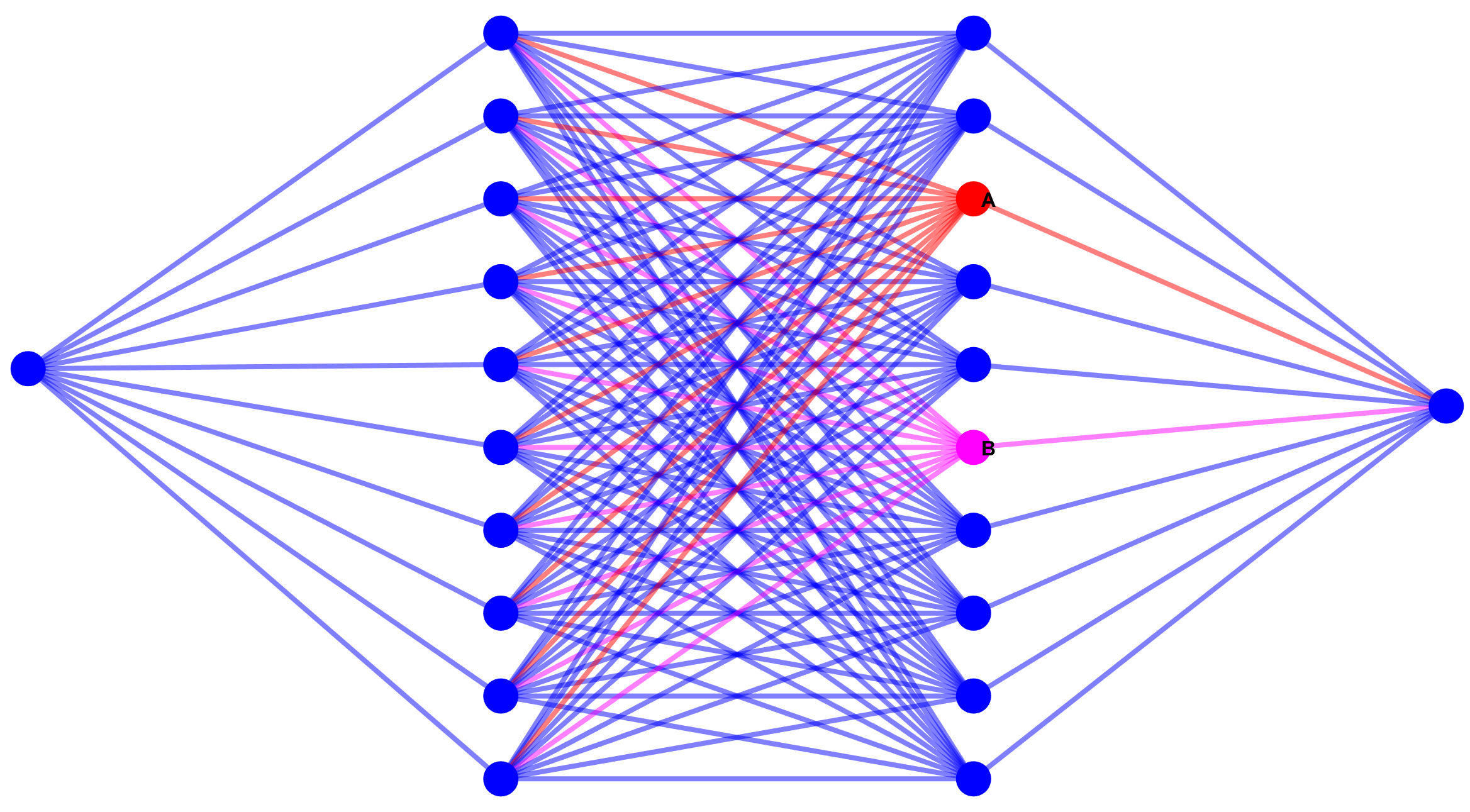}
\caption{A neural network with 2 hidden layers of width 10.  Switching the 
bias values and edges connected to node A (in red) with the those of 
node B (in magenta) results in a different ordering of the parameters, 
but the same neural network function.} \label{fig:netEx}
\end{figure} 
 
It is easy to see that there are many permutations of the parameters 
that result in the same network. This is illustrated in 
Figure \ref{fig:netEx} where a network is represented graphically with the 
nodes of the graph representing the biases and the edges of the graph 
representing the weights. Switching the bias values A (red) and B (magenta) 
results in the same network as long as the corresponding weights 
(also in red and magenta) are also switched.
This implies that it is possible to change the order of the parameters 
without changing the action of the function $\mathcal F$. In fact, 
for networks where each layer consists of $n_\ell$ neurons, there 
are $\sum_{\ell = 1}^{L-2} n_\ell!$ ways to rearrange the parameters 
and obtain the same network. 
For the relatively small network displayed in Figure \ref{fig:netEx} there 
are 7,257,600 different ways to permute the parameters and still obtain 
the same network. In the context of training, this means that for this 
particular example there are 7,257,600 different solutions to the optimization 
problem \eqref{eq:orgprob} that will result in the same neural network 
function $\mathcal F$.  This level of non-uniqueness is troubling and 
highly unsatisfactory. Apparently, this simple looking observation 
has not received any attention so far.

\medskip
\noindent
{\bf Problem formulation.} 
The parameter search space may be narrowed down by fixing the bias 
vector configurations to be monotonically increasing. This is achieved 
by augmenting the learning problem \eqref{eq:orgprob} 
with inequality constraints on the bias vectors in each layer. 
The resulting optimization problem that describes the training of the 
network is 
\begin{subequations} \label{eq:constrProb}
	\begin{alignat}{4}
		&\min_{\{W_\ell\}_{\ell=0}^{L-1}, \{b_\ell\}_{\ell=0}^{L-2}} && J(\{(y_L^i, S(u^i))\}_i,\{W_\ell\}_\ell,\{b_\ell\}_\ell) 
		\\
		 &\mbox{subject to } \quad && y_L^i = \mathcal F(u^i;(\{W_\ell\}, \{b_\ell\})) \, , \qquad   i = 1, \dots, N, \\
		  & && b_\ell^j \leq b_\ell^{j+1} \, , \qquad  j = 1, \dots, n_{\ell+1} - 1, \quad \ell = 0, \dots, L-2, \label{eq:ineqCon}
	\end{alignat}
\end{subequations}
where for a fixed $\ell$, the quantities $b_\ell^j$ are the entries of 
the bias vector $b_\ell$.  The added inequality constraints fix the 
configuration of the network, so that any permutation of the parameters 
(excluding the case when two adjacent bias values are equal) either 
violates the constraints or results in a different neural network.  

\medskip
\noindent
{\bf Outline of the paper.} 
Section~\ref{sec:RegProb} introduces a Moreau-Yosida regularization based 
algorithm to handle inequality constraints given in \eqref{eq:ineqCon}. 
This is followed by a convergence result of the Moreau-Yosida regularized 
problem to the original problem in Section~\ref{s:conv}. 
Section~\ref{s:resnets} describes the ResNets used for the examples
shown in Section~\ref{s:numerics}. In the first 
example the proposed algorithm is applied to a standard trigonometric 
function. This is followed by an application to a realistic application 
in chemically reacting flows, which are governed by stiff ordinary 
differential equations. 
The final example on the MNIST dataset within 
TensorFlow, illustrates that the presented approach can be incorporated 
in any of the existing machine learning libraries.

\section{Regularized problem} \label{sec:RegProb}

The inequality constraints in \eqref{eq:ineqCon} are difficult to 
implement directly. However, it is possible to implement them implicitly 
by augmenting the loss function $J$ with a penalty term similar to a 
Moreau-Yosida regularization. The well-known Moreau-Yosida regularization is 
frequently used to implement inequality constraints in the context of 
optimization problems with partial differential equations as constraints, 
see \cite{ItKu2008, NeTr2009, HiHi2009, AnBrVeWa2021}. Before introducing 
the loss function with bias order regularization, a more precise definition 
of the loss function is given. For training data $\{u^i, S(u^i)\}_{i=1}^N$, 
consider 
	\begin{equation}\label{eq:J}
		J: = \frac{1}{2N} \sum_{i=1}^N \|y_L^i - S(u^i)\|_2^2 + \frac{\lambda}{2}\sum_{\ell = 0}^{L-1}\big( \|W_\ell\|_1 + \|b_\ell\|_1 + \|W_\ell\|_2^2 + \|b_\ell\|_2^2\big),
	\end{equation}
where $\lambda>0$ is a regularization parameter. In order to fit $J$ into 
the framework introduced above, $b_{L-1}$ is taken to be the zero vector 
in $\mathbb R^{n_L}$.  
Even though a mean squared error term is used above to measure the 
approximation error of the neural network, this is easily generalizable to 
other terms such as cross-entropy, likelihood, etc. \cite{MR3617773}. 

Using the notation above to represent the entries of each bias vector, 
namely $b_\ell = (b_\ell^j)_{j=1}^{n_{\ell+1}} \in \mathbb R^{n_{\ell+1}}$, 
the new loss function is defined as 
\begin{equation}\label{eq:Jg}
	 J_\gamma := J + \frac{\gamma}{2} \sum_{\ell=0}^{L-2} \sum_{j = 1}^{n_{\ell+1}-1} \|\min \{ b_\ell^{j+1} - b_\ell^j, 0 \} \|_2^2, 
\end{equation}
where $\gamma$ is the so-called penalization parameter. The regularized 
optimization problem can now be written as 
\begin{subequations} \label{eq:regProb}
	\begin{alignat}{4}
		&\min_{\{W_\ell\}_{\ell=0}^{L-1}, \{b_\ell\}_{\ell=0}^{L-2}} J_\gamma(\{(y_L^i, S(u^i))\}_i,\{W_\ell\}_\ell,\{b_\ell\}_\ell) 
		\\
		 &\mbox{subject to } \quad y_L^i = \mathcal F(u^i; (\{W_\ell\}, \{b_\ell\})) \qquad i = 1, \dots, N.  
	\end{alignat}
\end{subequations}
Note that, even though it is not explicitly written in the formulation above, 
all of the variables $W_\ell$, $b_\ell$, and $y^i_L$ depend on the
parameter $\gamma$.

In Appendix \ref{app:1} the first order optimality conditions for this problem are 
derived where the DNN used is a Deep Residual Neural Net (ResNet).


\section{Convergence of $J_\gamma$}
\label{s:conv}

In order to show that as $\gamma \to \infty$, the 
minimum value of $J_\gamma$ converges to the minimum value of $J$ and the 
constraints \eqref{eq:ineqCon} are also satisfied, 
let $\theta$ represent the concatenation of all of the parameters which are 
being minimized, i.e. $\theta$ contains all of the entries of 
$\{W_\ell\}$ and $\{b_\ell\}$. Furthermore, assume that a fixed set of 
training data is being used and so the loss function $J$ given 
in \eqref{eq:J} can be rewritten as 
\[
J(\theta) := \frac{1}{2N} \sum_{i=1}^N \|\mathcal F(u^i; \theta) - S(u^i)\|_2^2 + \frac{\lambda}{2} \sum_{\ell = 0}^{L-1}\big( \|\theta\|_1 + \|\theta\|_2^2\big).
\]
Introducing the notation
\[
 g(\theta):=  \sum_{\ell=0}^{L-2} \sum_{j = 1}^{n_{\ell+1}-1} \|\min(b_\ell^{j+1} - b_\ell^j, 0)\|_2^2 \, ,
\]
the regularized loss function in \eqref{eq:Jg} may be rewritten as
\[
J_\gamma (\theta^\gamma) := J(\theta^\gamma) + \frac{\gamma}{2}g(\theta^\gamma).
\]
Now, the framework of \cite[Section 10.11]{Luenberger1969} can be used to 
show the following convergence results. The proof follows exactly as 
in \cite{Luenberger1969} after a transformation of notation. 

\begin{prop}\label{prop:conv}
Let $J_0$ be the minimum value attained from solving \eqref{eq:constrProb}. 
For each $\gamma$, let $\theta^\gamma$ be a minimizer of $J_{\gamma}$. 
The following hold
\begin{enumerate}
\item[(a)] $J_{\gamma} (\theta^\gamma) \geq J_{\tilde\gamma}(\theta^{\tilde\gamma})$, for 
	$\gamma \ge \tilde\gamma$; 
\item[(b)] $J_0 \geq J_{\gamma}(\theta^\gamma)$ for each $\gamma$;
\item[(c)] $\lim_{\gamma\to \infty} \frac{\gamma}{2} g(\theta^\gamma) = 0$. 
\end{enumerate}
\end{prop}

In particular, part $(c)$ shows that as $\gamma \to \infty$ the 
inequalities in \eqref{eq:ineqCon} are satisfied.

\begin{remark}
Typically, when a Moreau-Yosida regularization is implemented, a path-following technique is used to increase the size of $\gamma$ gradually.  This means that a sequence of optimization problems is solved for subsequently larger values of $\gamma$.   The initial $\gamma$ value is taken to be small, and the solution to the problem, $\theta^\gamma$, is used as the initial guess for the next optimization problem with a larger value of $\gamma$.   This path-following process continues until $\gamma$ is sufficiently large.   For the numerical examples presented below, a path-following technique was not used, and it was sufficient (as the results show) to solve each problem for a single value of $\gamma$. 
\end{remark}

\section{ResNets} 
\label{s:resnets}

As some of the examples shown below are obtained for Deep Residual Neural 
Networks (ResNets), a small description follows. 
Recall the definition of $\mathcal F$ in \eqref{eq:NNFunc}.  In the sequel, 
the layer functions will be denoted as 
$f_\ell = f_\ell( y_\ell; (W_\ell, b_\ell)): 
\mathbb R^{n_\ell} \to \mathbb R^{n_{\ell+1}}$, where the dependence of 
$f_\ell$ on $\sigma$ is not explicitly written.  With this representation the 
neural network can be viewed as an iterative progression of updating the 
output of each layer as
\[
y_{\ell + 1} = f_\ell(y_\ell; (W_\ell, b_\ell)) \qquad \ell = 0, \dots, L-1,
\]
with initial input $y_0$ and final output $y_L$.  As before, in order 
to preserve consistency, $b_{L-1}$ is taken to be the zero vector. 
One way to define the layer functions is 
\begin{alignat*}{5}
f_0(y_0) &:= \sigma(W_0 y_0 + b_0), \\
f_\ell(y_\ell) & := P_\ell y_\ell + \tau \sigma(W_\ell y_\ell + b_\ell) \qquad \ell = 1, \dots, L-2,\\
f_{L-1} &:= W_{L-1} y_{L-1},
\end{alignat*}
for matrices $P_\ell \in \mathbb R^{n_\ell \times n_{\ell+1}}$ and a 
scalar $\tau > 0$. 
With this definition and for $L>2$, $\mathcal F$ is termed as a Deep 
Residual Neural Network (ResNet).  If each $P_\ell$ is taken to be the 
identity matrix, which requires that the hidden layers have a uniform width, 
then this network can be viewed as a forward Euler discretization of an ODE.  
For more on these kinds of networks see 
\cite{HeZhReSu2016, RuHa2020, HaKhLoVe2020,antil2021novel,MR4060448}, among 
others. Notice that if $\tau = 1$ and $\{P_\ell\}_{\ell = 0}^{L-1}$ contains 
only zero entries, then $\mathcal F$ is a standard feedforward deep neural 
network \cite{MR3617773}.   

\section{Numerical Results}
\label{s:numerics}

In this section several examples are given that not only show the efficacy 
of this method, but also the advantages of using the method. In Section 
\ref{sec:sin}, a first example compares a single ResNet to learn the 
function $\sin(x)$ with and without bias ordering. This simple example 
shows that the method performs as desired, and in fact outperforms the 
same network trained without bias ordering. Following this, Section \ref{sec:chem}, a more complicated experiment is reported that 
use parallel ResNets to learn a model related to chemically reacting flows
\cite{brown2021novel}. This example shows that the proposed method is a 
useful technique for practical problems in machine learning.   Finally, in in Section \ref{sec:keras}, an example is described where the bias ordering regularization is applied to a classification problem using MNIST data and implemented in Keras.  This example shows that the proposed method is also suitable for Convolutional Neural Networks.

In Sections \ref{sec:sin} and \ref{sec:chem}
a BFGS optimization routine with 
Armijo line search is used during training to solve the optimization 
problem with or without bias ordering. In order to avoid overfitting, 
validation data is used with a patience of 400 iterations. This means 
that the training data is separated into two sets: training and validation 
data. The validation data is not used to update the weights and biases, 
rather it is used to measure the error of the network on unseen data 
(the validation data) during training.  If the validation error increases, 
the patience iterations provide a buffer during which training continues.  
If during these iterations, the validation error reaches a new minimum, 
then training continues as before, otherwise the training routine is 
terminated.
 
\subsection{ResNet to learn $\sin(x)$} \label{sec:sin}

For the first example a simple ResNet with 2 hidden layers of width 
50 is used to learn the function $\sin(x)$. The skip parameter $\tau$ is 
taken to be 1. For data, 1000 evenly spaced points from the interval 
$[0, 2\pi]$ are generated and then randomly split into training data 
(400 points), testing data (400 points), and validation data (200 points).

\begin{figure}[htb]
\center{
\includegraphics[width = \textwidth]{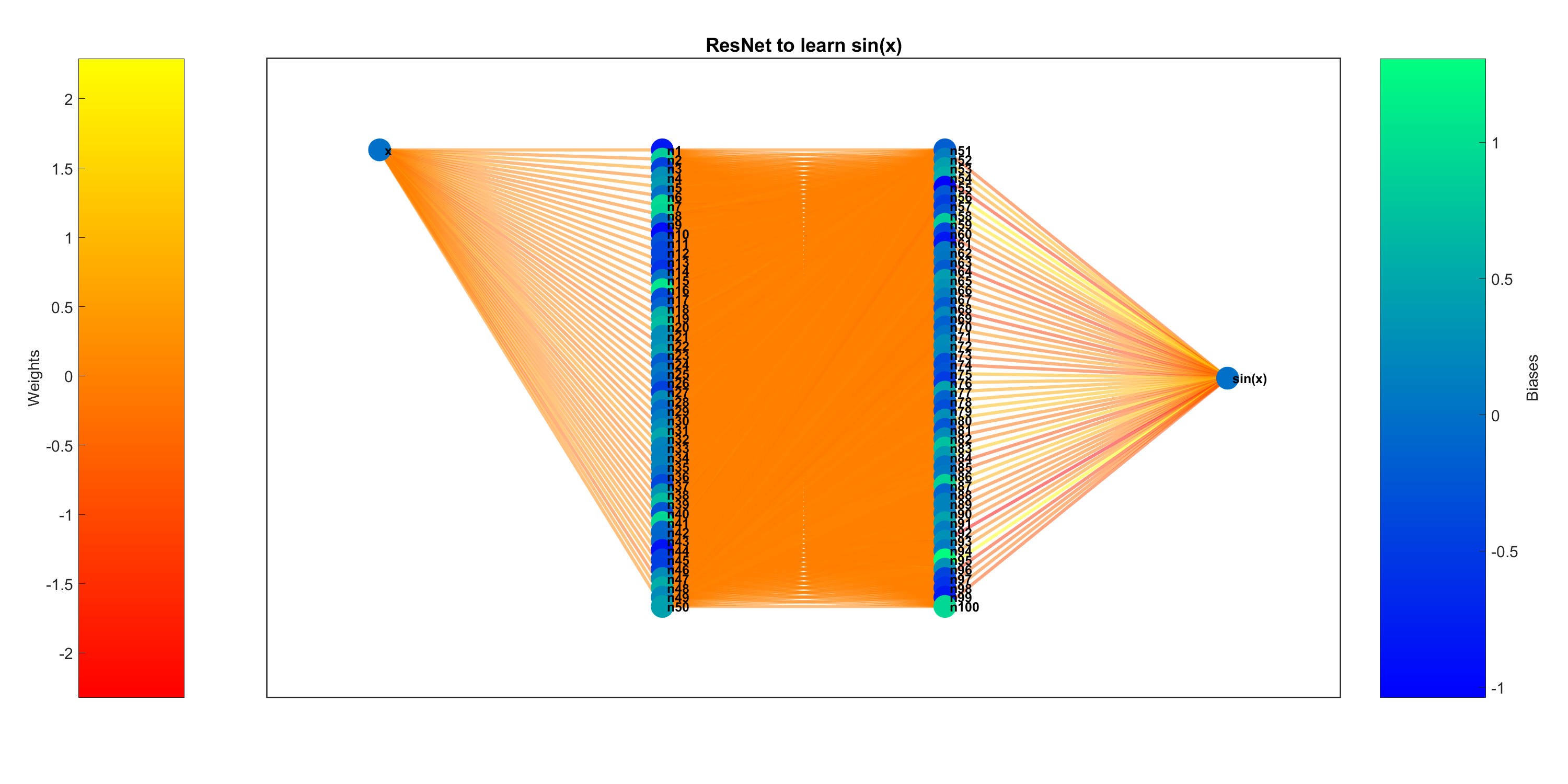}\\
\includegraphics[width = \textwidth]{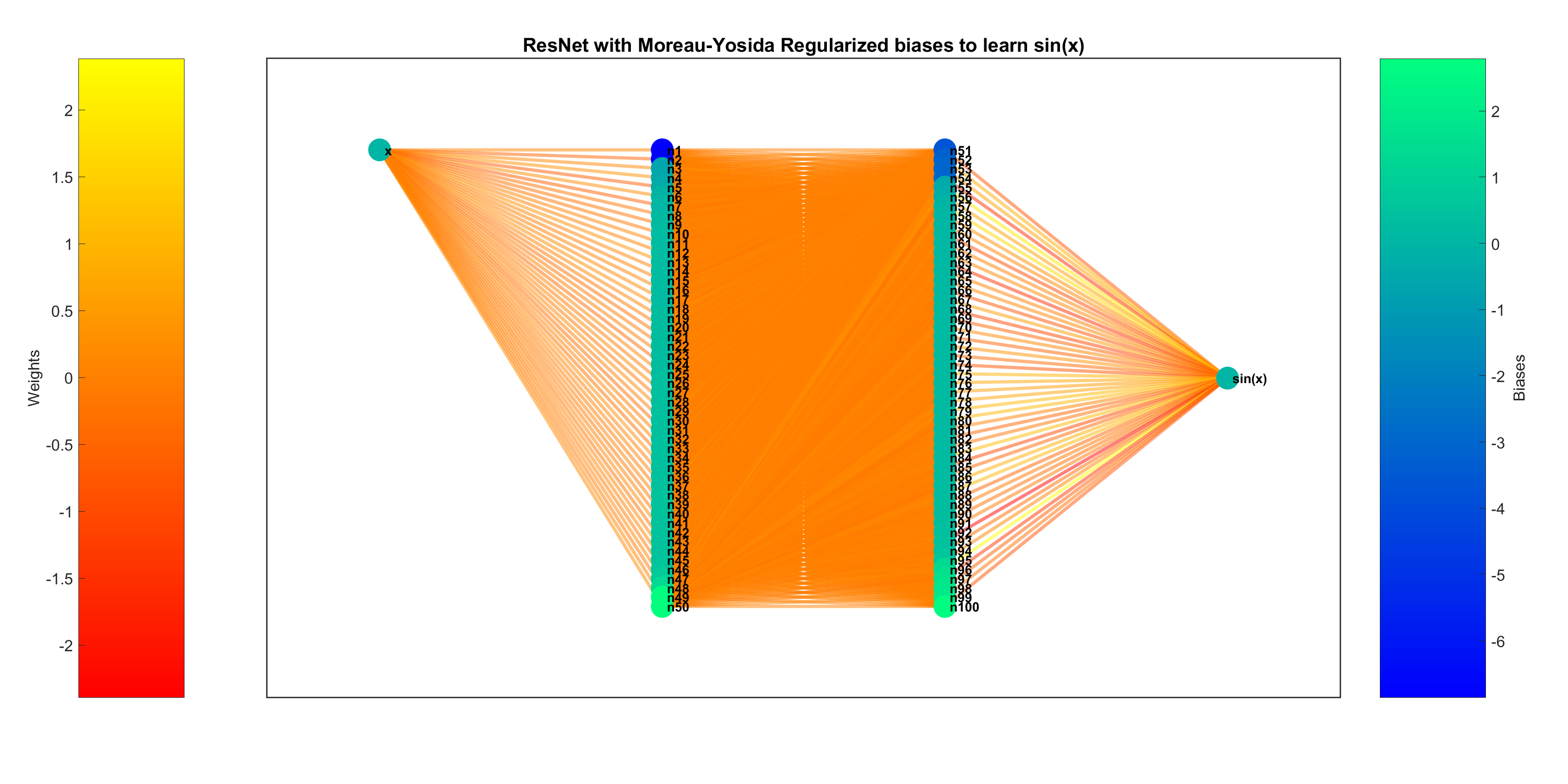}
\caption{ Network visualizations for the ResNets trained with loss function 
$J$ (top) and $J_\gamma$ (bottom). The network consists of 2 hidden layers 
with 50 neurons in each layer to learn $\sin(x)$. It is clear from the 
panels that the ordering of the  bias values is only enforced in the bottom 
panel.} \label{fig:1}
}
\end{figure}

The same experiment was performed twice, once using the loss function $J$, 
and once with the loss function $J_\gamma$ with $\gamma = 100$.  The 
resulting trained networks are shown in Figure \ref{fig:1}, where the 
bias values are represented by the neurons of the network, and the weights 
are represented by the connections between neurons.  The color for the input 
and output are set to zero.  In the case where the network was trained with 
loss function $J_\gamma$, the resulting biases were perfectly ordered and so 
the inequality constraints were satisfied.  This can also be seen in the 
bottom panel in Figure \ref{fig:1}. Note that for the network with 
unconstrained bias values, there are on the order of $10^{64}$ permutations 
of the parameters that will give the same ResNet approximation. For the 
ResNet with ordered bias values, however, any permutation of the parameters 
will either result in a different ResNet approximation, or violate the bias 
ordering.

\begin{figure}[htb]
\center{
\includegraphics[width =0.75 \textwidth]{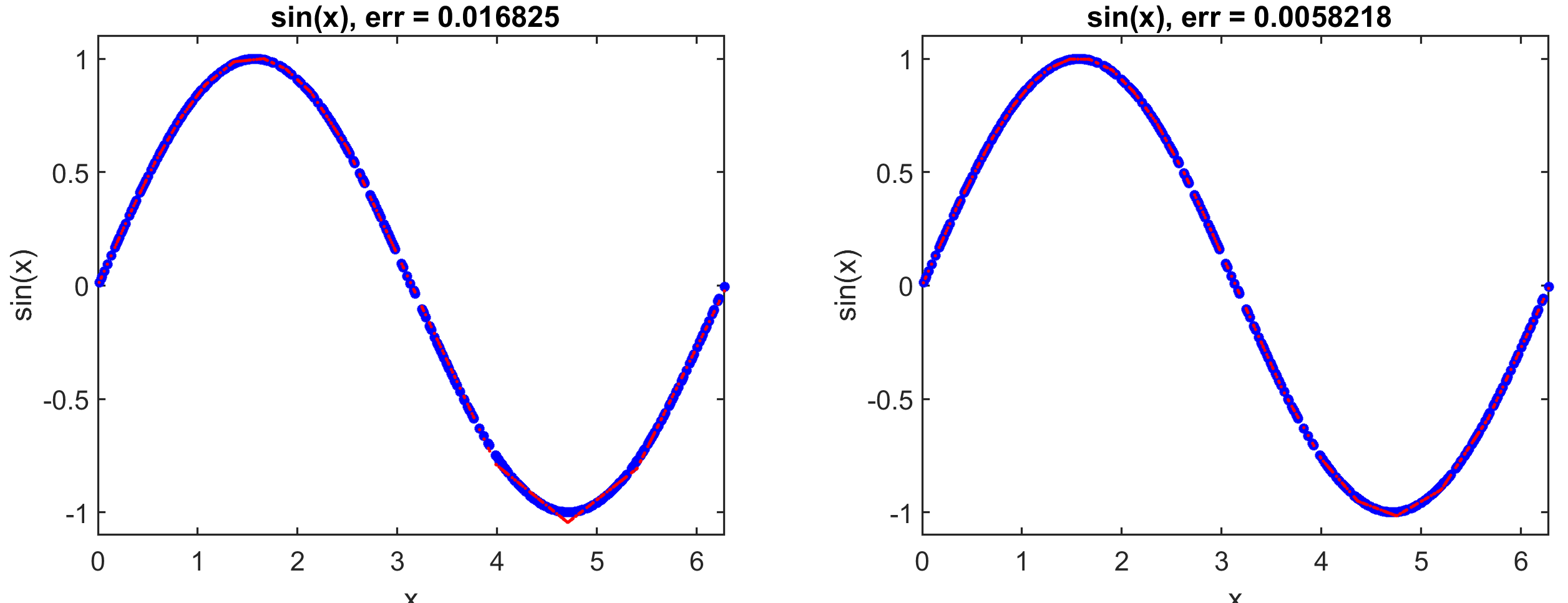}
\caption{Results from training ResNets  to learn the function $\sin(x)$.  
Exact values in blue, ResNet output in red. 
The left panel shows results for the standard ResNet with $J$ and the right 
panel shows results where bias ordering is enforced using $J_\gamma$. Both 
visual inspection and quantitative inspection of the error confirms that the 
proposed approach works better in this example.
} 
\label{fig:2}
}
\end{figure}

In Figure \ref{fig:2} the output of the two networks is compared on the 400 
test points, with exact values in blue and ResNet output in red. The left 
plot of Figure \ref{fig:2} shows the results for a network trained with a 
standard loss function $J$, while the plot on the right shows results from 
a network trained with the augmented loss function $J_\gamma$ which 
orders the biases. It is evident from the plots that the ResNet that implements 
bias order produced more accurate results. This is quantified above 
the plots with the relative error measured in the 2-norm. 
  
\subsection{Applications to chemically reacting flows} \label{sec:chem}

In this section experiments are presented involving a ResNet approximation of 
a system of stiff ODEs that model a reduced H$_2$-O$_2$ reaction. This reduced 
model (see \cite{PeHa1999}) tracks 8 species and temperature as they interact 
over time and is completely separated from any advection and diffusion in 
space.  For more information on this problem and more experiments using a 
parallel ResNet approximation see \cite{brown2021novel}.  

For the experiments included in this work, nine parallel ResNets with input 
dimension 10 and output dimension 1 (9 total output quantities) and 8 hidden 
layers of width 30 were trained on H$_2$-O$_2$ reaction data created by 
solving the stiff ODE system using CHEMKIN \cite{KeEtAl2000}.  Given an 
input vector representing the data at time $t_k$, each ResNet is trained to 
learn a single quantity (temperature, for example) at time $t_{k+1}$.  
More details can be found in \cite{brown2021novel}.  The generated data used 
to train and test these networks consists of thirteen subsets corresponding to 
initial conditions with a fixed equivalence ratio of 1 and different 
temperatures varying from 1200K to 2400K in increments of 100K.  
The parallel ResNets were trained on the data sets corresponding to initial 
temperatures 1200K, 1500K, 1800K, 2100K, and 2400K. 

The experiment above was performed twice, once with loss functions $J$  
(one for each parallel ResNet) as described in Section \ref{sec:RegProb}, and 
once with regularized loss functions $J_\gamma$ with $\gamma = 1000$.  All 
other network hyperparameters including the initial values of the weights and 
biases prior to training are kept the same for the two experiments.  
In Table \ref{tab:iter} the number of BFGS iterations used in training the 
parallel networks for both experiments described above is compared.  
The network with Moreau-Yosida regularization to order the bias values 
used fewer iterations to train in five of the nine networks.

\begin{table}[htb]
\begin{tabular}{|l|c|c|}
\hline
\multicolumn{3}{|c|}{Summary of BFGS iterations used during training}\\
\hline 
& with loss functions $J$ & with loss functions $J_\gamma$\\
\hline
Network 1 (temperature) & 2638 & 786 \\
Network 2 (O) & 2506 & 5125 \\
Network 3 (H) & 2772 & 1123 \\
Network 4 (OH) & 5957 & 4444\\
Network 5 (HO$_2$) & 230 & 540 \\
Network 6 (H$_2$O$_2$) & 886 & 1738 \\
Network 7 (H$_2$O) & 7267 & 583 \\
Network 8 (O$_2$) & 8872 & 2547\\
Network 9 (H$_2$) & 3719 & 4570\\
\hline
\end{tabular}
\caption{A comparision of the number of BFGS iterations used during training 
for the two sets of parallel ResNets.  The quantity that the ResNet is 
learning (output) is written in parentheses.} \label{tab:iter}
\end{table}

Recall from Proposition~\ref{prop:conv}, that the Moreau-Yosida regularization 
approach will be enforcing the ordering \eqref{eq:ineqCon} approximately.
  
The loss functions $J_\gamma$ penalize the violation of the ordering, but do 
not strictly implement the order itself.  Even so, the biases in eight of the 
nine parallel networks were ordered perfectly for the networks that trained 
with the loss functions $J_\gamma$.   The only violation of the bias ordering 
occured in Network 9, the network used to learn H$_2$.  In the first hidden 
layer of this network neurons six and seven were ordered incorrectly.  In the 
trained network, the value of neuron 6 was approximately -0.0175652 and the 
value of neuron 7 was approximately  -0.0175660, and therefore the size of 
the order violation was -8e-7, which is negligible.  To reiterate, in each 
of the 9 networks, there were 240 bias values (30 per hidden layer) for a 
total of 2,160 different biases.  In all of these biases only a single pair 
(negligibly) violated the monotonic ordering using the proposed method.  

\begin{figure}[htb]
\center{
\includegraphics[width =0.45 \textwidth]{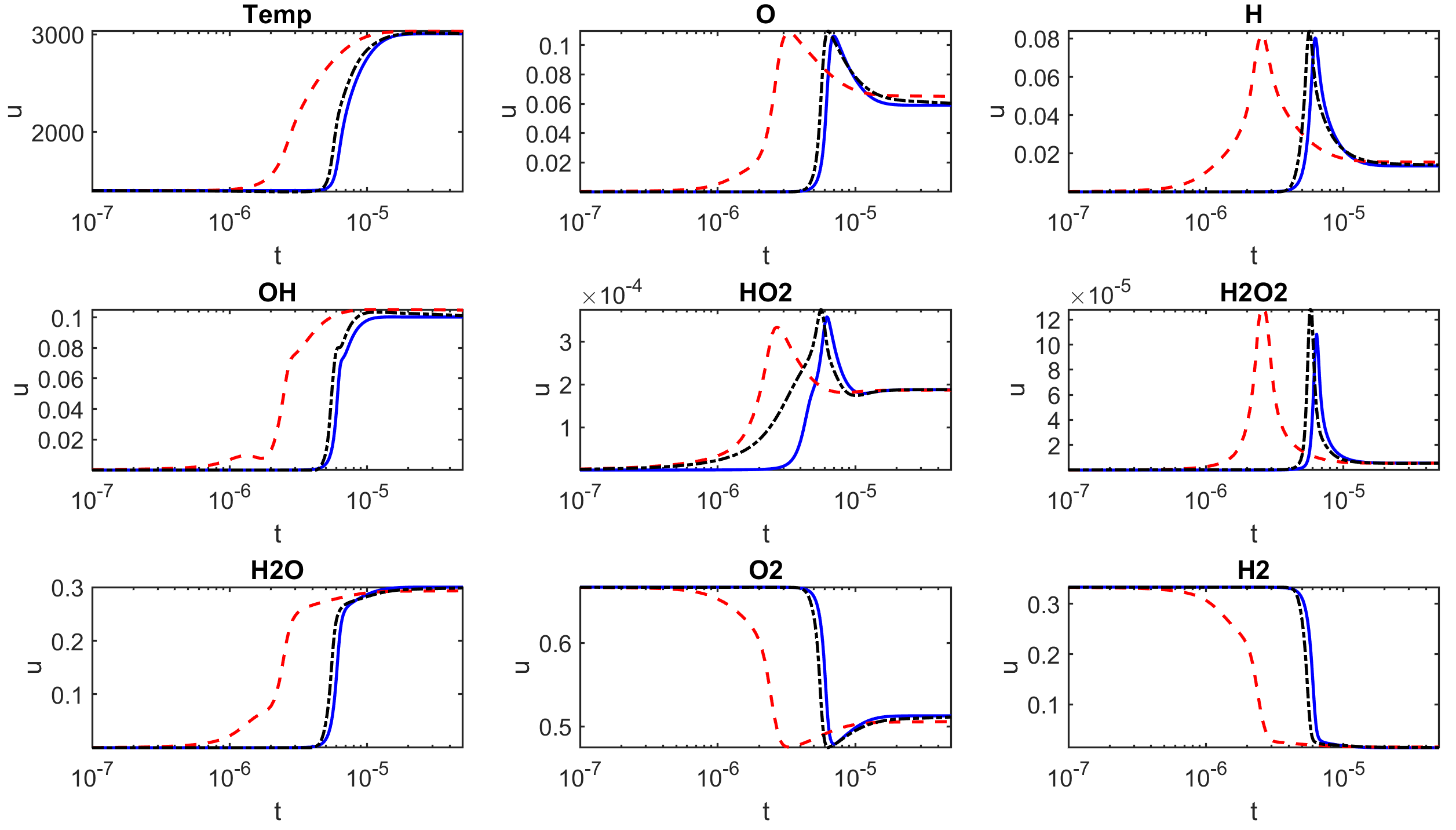} \quad
\includegraphics[width = 0.45\textwidth]{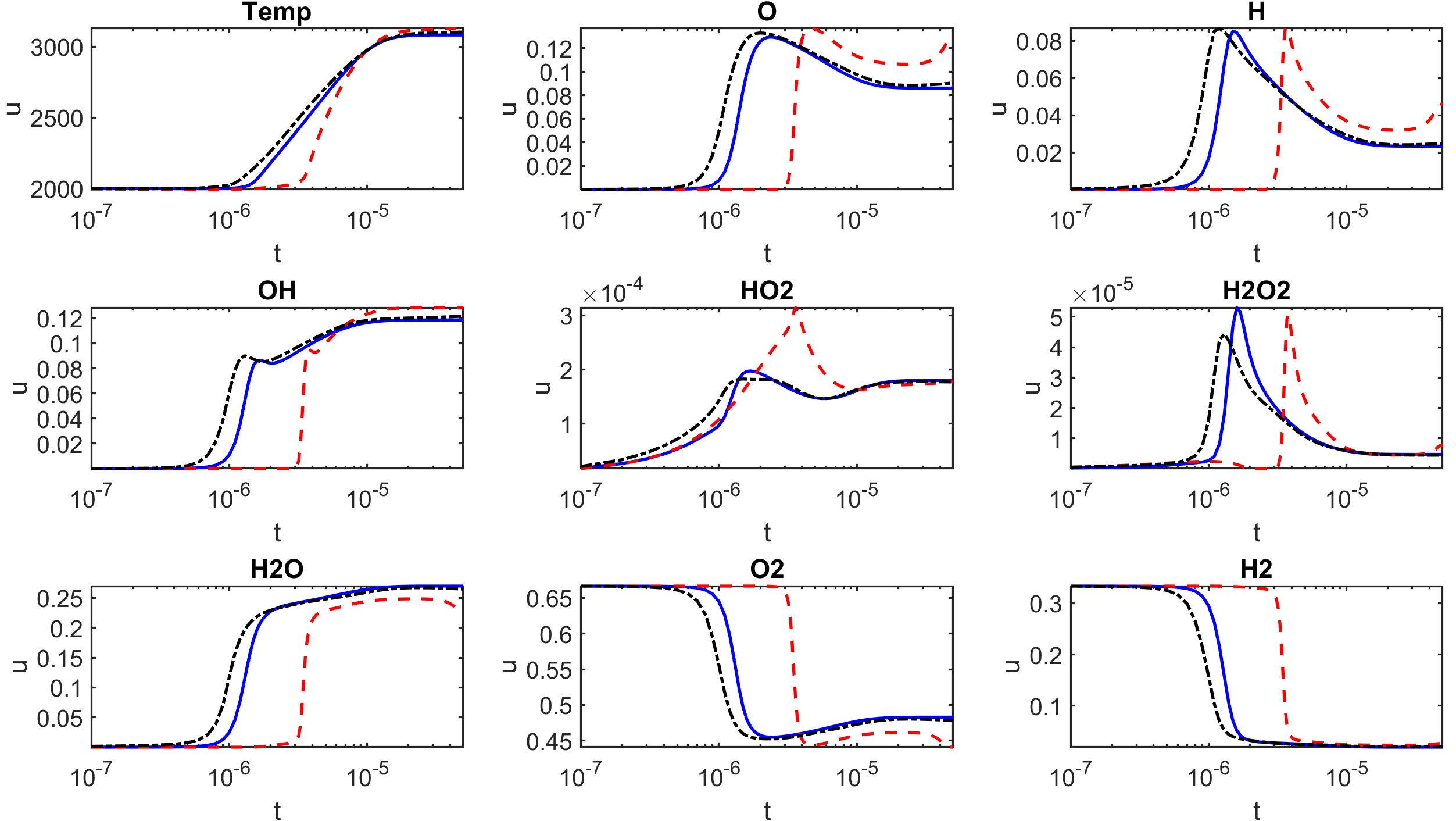}
\caption{Results from training parallel ResNets on a reduced H$_2$-O$_2$ 
reaction model.  Values from known data in blue, results from ResNets trained 
with traditional loss functions in red, and results from ResNets trained 
with loss functions with a Moreau-Yosida penalization term to order the biases 
in black. Clearly, the proposed  approach outperforms the existing one.}
\label{fig:3}
}
\end{figure}

Figure \ref{fig:3} compares the results of the two sets of parallel ResNets 
tested with initial conditions with initial temperature 1400K (left set of 
plots) and 2000K (right set of plots).   Note that the networks were not 
trained on this data.  To test the networks only the initial condition comes
from the CHEMKIN data.  The output of the ResNets from the initial condition 
are then combined and used as the input to the parallel ResNets for the next 
timestep.  This process is repeated for the duration of the reaction.  
To compare the results, the known CHEMKIN data are represented in blue, the 
results from the networks trained with $J$ are represented with dashed red 
lines and the results for the networks trained with $J_\gamma$ are 
represented with dash/dotted black lines. It can clearly be seen that the 
results of the parallel ResNets trained with $J_\gamma$ match the data 
more closely.  Furthermore, note that the $x-$axis (time axis) is scaled 
logarithmically in order to display the details of H$_2$O$_2$ which happen 
quickly and early in the reaction.  Therefore, while the results in black 
anticipate the reaction, they only do so slightly, on the order 
of $10^{-6}$ seconds.

\subsection{A Convolutional Neural Network classification problem} \label{sec:keras}

For this example, a convolutional neural network is constructed using Keras to solve a classification problem using the MNIST data set.  The purpose of this example is to show that the flexibility of the bias ordering technique and that it can be incorporated into existing neural network software.   Unlike the other examples presented, the neural network is not a ResNet.  Instead the network consists of 4 convolutional layers followed by two dense layers.  The diagram in Figure \ref{fig:netEx} can also be used to describe a convolutional layer, where the nodes of the graph still represent the bias value, but now the edges represent a convolutional filter or convolutional kernel.  In this way, the bias ordering technique easily applies to these networks as well.  

The network used for this example consists of four convolutional layers followed by two dense layers.  The first two convolutional layers consist of 32 convolutional filters of width 3.  Convolutional layers three and four both consist of 64 convolutional filters of width 3.  The first dense layer consists of 512 neurons, and the second dense layer (the output layer) consists of 10 neurons since there are 10 digits (or classes) in the MNIST data.  Unlike the other examples, bias values are included on the output layer of this network.  A softmax activation function is used for the output layer, while a ReLU activation function is used for the other five layers.    Max pooling is used after layers 2 and 4.  In order to avoid overfitting, dropout is used in the network.   

Two versions of the network described above are constructed using Keras with a Tensorflow backend.   Both of the networks use a categorical cross entropy loss function, and the difference between the two is that one of the networks also implements the Moreau-Yosida regularization term through a custom bias regularizer.  The regularizer is defined using the following lines in Python
\begin{lstlisting}[language=Python]
def MY_regularizer(bias, gamma=100):
    bias_length = tf.size(bias)
    bias_diff = bias[1:bias_length]-bias[0:bias_length-1]
    bias_min = tf.math.minimum(bias_diff,0)
    bias_norm = tf.math.reduce_sum(tf.math.square(bias_min))
    return 0.5*gamma*bias_norm 
\end{lstlisting}
In the above code , the value of the regularization parameter is taken to be 100, but this is easily customized.  Once the regularizer is defined, it can be implemented into any layer with the {\tt bias\_regularizer} option.   

\begin{figure}[htb]
\center{
\includegraphics[width=0.45 \textwidth]{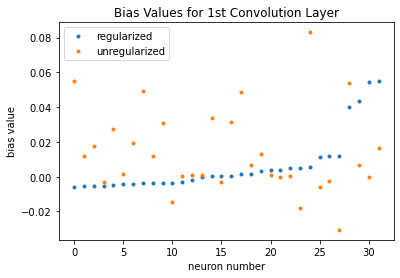} \quad
\includegraphics[width=0.45 \textwidth]{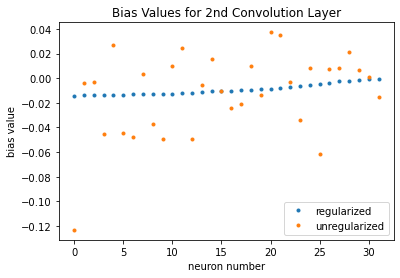} \\
\includegraphics[width=0.45 \textwidth]{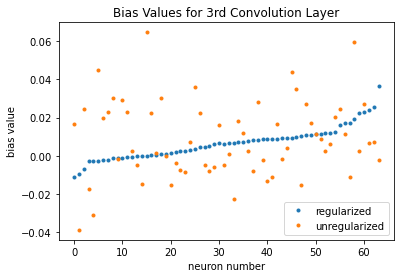} \quad
\includegraphics[width=0.45 \textwidth]{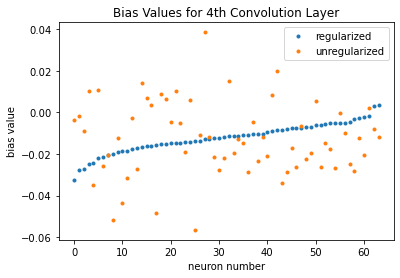} \\
\includegraphics[width=0.45 \textwidth]{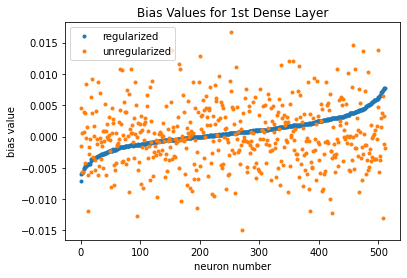} \quad
\includegraphics[width=0.45 \textwidth]{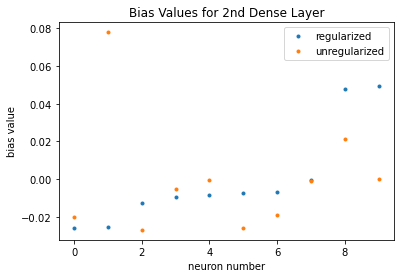} 
\caption{The bias values by layer for the two neural networks described in Section \ref{sec:keras}.  These plots show the successful implementation of the Moreau-Yosida bias order regularization in Keras with a Tensorflow backend.} \label{fig:layers}
}
\end{figure}

Both of the networks were trained on 60,000 samples from the MNIST data set.  The loss functions were minimized using stochastic gradient descent with a learning rate of 0.01.  Batch normalization was used during training with a batch size of 32.  Each network trained for 5 epochs before the networks were tested on 10,000 MNIST samples. In Figure \ref{fig:layers} the bias values for each layer of the two networks are plotted.  It can clearly be seen that the regularization resulted in successfully ordering the bias values in each layer.  The accuracy of the two trained networks were comparable.  For the networks corresponding to the plots in Figure \ref{fig:layers}, the network with the bias regularization had an accuracy of 98.7\% and the network without regularization had an accuracy of 98.66\%. 

\section{Conclusions}

A method to reduce the very large search space of equivalently optimal
neural nets has been introduced. The key idea is to try to enforce 
the biases to be monotonically 
increasing in each layer of neurons. This is accomplished
by a Moreau-Yosida regularization-based approach to solve the resulting 
optimization problem. The convergence of the regularized problem has been
proven. The benefit of this new approach has been demonstrated not only
on simple approximation cases but also on a realistic problem arising 
in chemically reacting flows. The numerical experiments presented also 
show that the regularization technique is effective for relatively 
small values of the penalization parameter $\gamma$. 

\appendix
\section{Derivation of the first order optimality conditions} \label{app:1}

In the context of constrained optimization, the problem \eqref{eq:regProb} is typically solved by using a gradient based method.  Indeed this is the approach taken in the numerical experiments.  Each evaluation of the gradient requires solving the state equation or forward problem and solving the adjoint equation problem.  In what follows, the state and adjoint equations will be derived for the problem in \eqref{eq:regProb} where the neural network being used is a Deep Residual Neural Network (ResNet), as introduced in Section \ref{s:resnets}.  

The state and adjoint equations as well as the gradient will be derived by using the Lagrangian approach.  For brevity, the following will be written for a single input $u$, rather than for a set of training data.  Similar to Section \ref{s:conv}, $\theta$ is used to represent the concatenation of the weights and bias, i.e., the parameters being optimized.  For appropriate adjoint variables $\psi = (\psi_j)_{j=1}^L$ the Lagrangian functional corresponding to \eqref{eq:regProb} is 
\begin{multline*}
		\mathcal{L}(u,\theta,\psi) := J_\gamma(\theta) +	\langle y_1 -\tau\sigma(W_0 u +  b_0), \psi_1  \rangle \\
		+ \sum_{j=2}^{L-1} \Bigg\langle y_{j} - P_{j-1}y_{j-1} -  
		\tau\sigma(W_{j-1}y_{j-1} + b_{j-1}), \psi_j  \Bigg\rangle 
		+ \langle y_L - W_{L-1} y_{L-1}, \psi_L  \rangle	,
	\end{multline*}

From here, the state and adjoint equations result
from evaluating the derivatives of the Lagrangian with respect to $y_j$ and $\psi_j$ at a stationary point.  Furthermore, the gradient is derived by taking the derivatives of the Lagrangian with respect to $\theta$.  
	\begin{subequations}\label{eq:oc_fDNN}
		\begin{enumerate}[(i)]
			\item State Equation. 
				\begin{align} \label{eq:oc_state}	
					\begin{aligned}
						y_1 &= \tau \sigma(W_0 u + b_0),  \\
						y_{j} &= P_{j-1}y_{j-1}  + \tau \sigma(W_{j-1}y_{j-1} + b_{j-1}),\qquad   2\leq j\leq L-1 , \\
						y_L   &= W_{L-1} y_{L-1} . 
					\end{aligned}
				\end{align}			
			\item Adjoint Equation.			
				\begin{equation}\label{eq:oc_adjoint}
					\begin{aligned}
						\psi_j &= P_{j}^T\psi_{j+1} -  \tau\left[-W_j^{T}\left(\psi_{j+1} \odot 	\sigma^{\prime}\left(W_{j} y_{j+1} + b_j \right)\right)\right] \qquad j=L-2,\hdots,1\\
						\psi_{L-1}&=-W_{L-1}^T \psi_L, \\
						\psi_L &= -\partial_{y_L} J_\gamma(\theta).
					\end{aligned}
				\end{equation}				
			\item Derivative with respect to $\theta$.
				\begin{equation}\label{eq:design}
					\begin{aligned}
						\partial_{W_{L-1}}\mathcal{L} =&-\psi_L\;y_{L-1}^T+\partial_{W_{L-1}}J_\gamma(\theta) \\ = &\partial_{y_L}J_\gamma(\theta)\;y_{L-1}^T+\partial_{W_{L-1}}J_\gamma(\theta),\\
						\partial_{W_{j}}\mathcal{L}= &-y_{j}\:\left(\psi_{j+1} \odot \sigma^{\prime}(W_{j} y_{j} + b_j)\right)^{T}+\partial_{W_j}J_\gamma(\theta)  \qquad&& j = 0,...,L-2,  \\
						\partial_{b_{j}}\mathcal{L}= &-\; \psi_{j+1}^T\;\sigma^{\prime}(W_{j} y_{j} + b_{j})+\partial_{b_j}J_\gamma(\theta) \qquad &&j = 0,...,L-2 \, .
					\end{aligned}
				\end{equation}
		\end{enumerate}			
	\end{subequations}

In the jargon of machine learning the equations in \eqref{eq:oc_adjoint} are called back propagation.  The gradient is represented by the right hand side of \eqref{eq:design}, where the contributions from the Moreau-Yosida regularization terms enter in the term $\partial_{b_j} J_\gamma (\theta)$.

\bibliographystyle{plain}
\bibliography{references}

\end{document}